\documentclass[a4paper,12pt]{article}

\usepackage{amsmath,amssymb,amsthm}
\usepackage{amsfonts}
\usepackage{multicol}
\usepackage{color}

\newtheorem{theorem}{Theorem}[section]
\newtheorem{lemma}[theorem]{Lemma}
\newtheorem{e-proposition}[theorem]{Proposition}
\newtheorem{corollary}[theorem]{Corollary}
\newtheorem{e-definition}[theorem]{Definition\rm}

\newcommand{\charI}{1 \hspace*{-1mm} {\rm l}}
 
\def\ignore#1{}
\def\eq{\begin{equation}}
\def\en{\end{equation}}
\def\eqa{\begin{eqnarray}}
\def\ena{\end{eqnarray}}
\def\eqs{\begin{eqnarray*}}
\def\ens{\end{eqnarray*}}
\def\bZ{{\mathbb Z}}

\def\re{{\mathbb R}}

\def\non{\nonumber}

\def\tg{{\tilde g}}
\def\s{\sigma}
\def\l{\lambda}
\def\d{\delta}
\def\Ref#1{(\ref{#1})}

\def\Eq{\ =\ }
\def\Le{\ \le\ }

\def\h{\eta}
\def\t{\tau}

\def\e{\varepsilon}

\def\a{\alpha}

\def\gg{{\mathcal G}}

\def\pr{{\mathbb P}}
\def\integ{{\mathbb Z}}
\def\ex{{\mathbb E}}
\def\tg{{\tilde g}}
\def\tgr{\tg_{nv_c,r}}
\def\nvfl{{\lfloor nv_c \rfloor}}
\def\ncfl{{\lfloor nc \rfloor}}
\def\nvan{{\langle nv_c \rangle}}
\def\siz{\sum_{i\in\integ}}
\def\sprz{{\sup_{r\in\integ}}}

\def\btd{\bigtriangledown}
\def\sln{\sqrt{\log n}}
\def\ui{^{(1)}}
\def\ut{^{(2)}}
\def\uh{^{(3)}}
\def\snl{{\lfloor \sqrt n \rfloor}}
\def\snu{{\lceil \sqrt n \rceil}}
\def\siii{\sum_{i=I_1}^{I_2-1}}
\def\ii{{(I_1,I_2)}}
\def\dtv{d_{TV}}
\def\aa{{\mathcal A}}
\def\law{{\mathcal L}}
\def\L{\Lambda}

\def\D{\Delta}
\def\m{\mu}
\def\giv{\,|\,}
\def\Blb{\left\{}
\def\Brb{\right\}}
\def\tS{{\widetilde S}}
\def\Po{{\rm Po\,}}
\def\bone{{\mathbf 1}}
\def\nin{\noindent}
\def\half{{\textstyle\frac12}}
\def\sjz{\sum_{j\in\integ}}
\def\Bl{\Bigl(}
\def\Br{\Bigr)}
\def\Giv{\,\Big|\,}
\def\brak#1{{\rm (#1)}}

\begin{document}

\title{Local limit approximations for \\ Markov
population processes}

\author{Sanda N.\ Socoll~and A.\ D.\ Barbour\footnote{Angewandte Mathematik, Universit\"at Z\"urich, 
Winterthurertrasse 190, CH-8057 Z\"URICH; 
work supported in part by Schweizerischer Nationalfonds Projekte Nr.\ 
20--107935/1 and 20--117625/1.}\  
\\ 
Universit\"at Z\"urich}

\date{}
\maketitle{}

\begin{abstract}
The paper is concerned with the equilibrium
distribution~$\Pi_n$ of the $n$-th element in a sequence of continuous-time 
density dependent Markov processes on the integers.
Under a $(2+\a)$-th moment condition on the jump distributions, 
we establish a bound of order $O( n^{-(\a+1)/2}\sqrt{\log n})$ 
on the difference between the point probabilities of~$\Pi_n$ and those of
a translated Poisson distribution with the same variance. Except for the
factor $\sln$, the result is as good as could be obtained in the simpler
setting of sums of independent integer-valued random variables.
Our arguments are based on the Stein-Chen method and coupling.
\end{abstract} 


{\it Keywords:}\ {continuous-time Markov jump process; equilibrium distribution; 
\hfil\break\hglue2.58cm point probabilities;
 Stein--Chen method; coupling}


\setcounter{equation}{0}
\section{Introduction}

Density dependent Markov population processes, in which the transition rates depend
on the density of individuals in the population, have proved widely useful as models 
in the social and life
sciences: see, for example, the monograph of Kurtz~(1981), in which approximations in terms of
diffusions are extensively discussed, in the limit as the typical population size~$n$ tends to
infinity. In the present paper we consider local approximation to their equilibrium 
distributions~$\Pi_n$. In Socoll~\&~Barbour~(2008) [SB], total variation approximation
to~$\Pi_n$ by a suitably translated Poisson distribution was shown to be accurate to
order $O(n^{-\a/2})$, provided that the jump distributions satisfy a $(2+\a)$-th moment condition
for some $0 < \a \le 1$. Here, we examine the approximation of point probabilities by
those of the same translated Poisson distribution, and show 
in Theorem~\ref{limth} that, under the same
assumptions, the error is now of order $O(n^{-(\a+1)/2}\sqrt{\log n})$. This is only
worse by the logarithmic factor than the best that can be obtained under comparable
conditions for sums of independent integer valued random variables.

A key ingredient in the proof of total variation approximation in~[SB] was to show
that the total variation distance between~$\Pi_n$ and its unit translate $\Pi_n*\d_1$
is of order~$O(n^{-1/2})$.  Here, we need to establish a local limit analogue of 
this theorem.  We prove in Section~2 that the differences between the point 
probabilities of~$\Pi_n$ and those of its unit translate are uniformly bounded by a
quantity of order $O(n^{-1}\sqrt{\log n})$. An important step in proving this is to
establish that, for some $U\ge1$, the difference between $\pr[Z_n(t) = k+1 \giv Z_n(0)=i]$
and $\pr[Z_n(t) = k \giv Z_n(0)=i-1]$ is of order  $O(n^{-1}\sqrt{\log n})$, uniformly
for~$i$ in a set~$I$ such that $\Pi_n(I^c) = O(n^{-1})$.  This is achieved by a
pathwise comparison of probability densities, and using a martingale concentration
inequality.  Note that, for sums of independent random variables, the corresponding
difference is always zero, so that this problem does not arise there.

The proof of Theorem~\ref{limth} is undertaken in Section~3.  The argument relies on
the Stein--Chen method (Chen, 1975) and Dynkin's formula, exploiting the particularly nice 
properties of the solutions to the Stein--Chen equation for one point subsets of~$\mathbb Z_+$. 

\subsection*{Preliminaries}\label{prelims}

For each $n\in\mathbb N$, let $Z_{n}(t)$, $t\ge0$,  be an 
irreducible continuous time pure jump Markov process taking values in $\mathbb Z$, with transition
rates given by 
$$
  i\ \to\  i+j \quad \mbox{ at rate }\quad n\lambda_j\Big(\frac{i}{n}\Big),\qquad i \in {\mathbb Z},\ 
  j\in \mathbb Z\setminus \{0\},
$$ 
where the $\lambda_j(\cdot)$ are prescribed functions on $\mathbb R$.
We then define the `overall jump rate' of the process $n^{-1}Z_{n}$ at 
$z \in n^{-1}{\mathbb Z}$ by
$$
  \L(z) \ :=\ \sum_{j\in \mathbb Z\setminus \{0\}} {{\lambda}_{j}}(z),
$$
its `average growth rate' by
$$
  F(z)\ :=\ \sum_{j\in \mathbb Z\setminus \{0\}}j{{\lambda}_{j}}(z),
$$ 
and its `quadratic variation' function by $n^{-1}{\sigma}^2(z)$, where 
$$
  {\sigma}^2(z)\ =\ \sum_{j\in \mathbb Z\setminus \{0\}} j^2  {{\lambda}_{j}}(z),
$$
assumed to be finite for all~$z\in\re$.

We make the following assumptions on the functions ${\lambda}_{j}$; they
are discussed at greater length in~[SB].
\\[1ex] 
{\bf \small A1:} There exists a unique $c$ satisfying $F(c)=0$; furthermore, $F'(c) < 0$
    and, for any $\eta>0$, $\mu_{\eta}:= \inf_{|z-c|\geq {\eta}}|F(z)|>0$. 
\\[1ex] 
{\bf \small A2:} 
   For each $j\in {\mathbb Z}\setminus \{0,\}$, there exists $c_j\ge0$ such that 
    \begin{equation} \label{lamdaj}
      \lambda_j(z)\leq c_j(1+|z-c|), \qquad z\in \mathbb R,
    \end{equation} 
    where the~$c_j$ are such that, for some $0 < \a \le 1$,
    $$
       \sum_{j \in \mathbb Z \setminus \{0\}}|j|^{2+\alpha}c_{j} =: s_\a <\infty.
    $$
\ignore{
The moment condition on the~$c_j$ in Assumption~A2\,(a) plays the same r\^ole as
the analogous moment condition in the Lyapounov central limit theorem.  Under
this assumption, the ideal rate of convergence in the usual central limit
approximation is the rate~$O(n^{-\a/2})$ that we had established for our total
variation approximation.
Assumption A2\,(b) is important for establishing the smoothness of the equilibrium
distribution~$\Pi_n$.  If, for instance, all jump sizes were multiples of~$2$, the
approximation that we are concerned with would not be accurate in total variation.
}
{\bf \small A3:}  
  There exist $\varepsilon>0$ and $0 < \delta \le 1$ and a set $J \subset 
     {\mathbb Z}\setminus \{0\}$ with $1 \in J$ such that
     \eqs
       &&\inf_{|z-c|\leq \delta}\lambda_j(z)\ \geq\ \varepsilon \lambda_j(c) \ > 0,\ \ \ j\in J;\\
       &&\lambda_j(z)=0 \  \mbox{for all}\   |z-c|\leq \delta,\ \ \ j\notin J.
     \ens   
\ignore{ 
Assumptions A2\,(a) and A3 imply in particular that the series  
$\sum_{j \in \mathbb Z \setminus \{0\}}j\lambda_j(z)$ and 
$\sum_{j \in \mathbb Z \setminus \{0\}}j^2\lambda_j(z)$  are uniformly convergent 
on~$|z-c| \le \d$,  and that their sums, 
$F$ and ${\sigma}^2$ respectively, are continuous there. They also imply that
$$
  \sum_{j\in \mathbb Z\setminus \{0\}}|j|n \lambda_j(i/n)
     \Eq O(|i|),\ \ \ |i|\to \infty,
$$
so that the process $Z_n$ is a.s.\ non-explosive, in view of 
Hamza \& Klebaner~(1995, Corollary~2.1).
}
\\[0.1ex] 
{\bf \small A4:} 
  \ (a)\ For each $ j\in J$, ${\lambda}_{j}$  is of class $C^2$ on~$|z-c| \le \d$.\\[0.3ex]
\mbox{} \qquad  (b)\  For $\delta$ as in A3,  
$$
   L_1 \ :=\ 
    \sup_{j\in J}\frac{\|{{\lambda}_j^{\prime}}\|_{\delta}}{\lambda_j(c)}
      \ <\ \infty;\qquad 
    L_2\ :=\ 
   \sup_{j\in J} \frac{\|{{\lambda }_j^{\prime \prime}}\|_{\delta}}{|j|\lambda_j(c)}
        \ <\ \infty,
$$
where   $\|f\|_{\delta}:=\sup_{|z-c|\leq \delta}|f(z)|$. \\[0.4ex] 
\ignore{
This assumption implies in particular, in view of Assumptions A2--A3, that the series 
$\sum_{j \in \mathbb Z \setminus \{0\}}j\lambda^{\prime}_j(z)$ and 
$\sum_{j \in \mathbb Z \setminus \{0\}}j^2\lambda^{\prime}_j(z)$ are uniformly 
convergent on $|z-c|\leq \delta$,
that their sums are 
$F^{\prime}$ and~$(\s^2)^{\prime}$  respectively, and that
$F$ and $\s^2$  are of class $C^1$ on $|z-c|\leq \delta$.
}
\ignore{
This assumption implies, in view of A2--A3, that the series 
$\sum_{j \in \mathbb Z \setminus \{0\}}j\lambda^{\prime \prime}_j(z)$ 
is uniformly convergent on $|z-c|\leq \delta$, its sum is
$F^{\prime \prime}$, 
and~$F$ is of class $C^2$ on $|z-c|\leq \delta.$}

In [SB], it is shown that, under assumptions A1--A4,
\eq\label{dtv-thm}
  d_{TV}(\widehat{\Pi}_n, \widehat{\rm Po}(nv_c)) \Eq O(n^{-\a/2}), 
\en 
where $v_c := \s^2(c)/\{-2F'(c)\}$, $\widehat{\rm Po}(nv_c)$ denotes the centred Poisson distribution 
$$
  \widehat{\rm Po}(nv_c)\ :=\ {\rm Po}(nv_c)*\d_{-\lfloor nv_c\rfloor},
$$ 
and $\widehat{\Pi}_n$ denotes the centred
equilibrium distribution ${\Pi_n}*\d_{-\lfloor nc\rfloor}$; 
here, $\d_r$ denotes the point mass at~$r$, and~$*$ denotes convolution. 

In this paper we prove the following complementary local limit approximation.

\begin{theorem} \label{limth} Under Assumptions A1--A4, there exists a 
constant $C>0$ such that 
 $$
  \sup_{k\in \mathbb Z} \Bigl|\widehat{\Pi}_{n}(k)-\widehat{\rm Po}(nv_c)\{k\}\Bigr|\Le
     C n^{-(\a+1)/2}\sqrt{\log n}.
$$ 
\end{theorem}
\nin This theorem shows that, even at the level of point probabilities, the approximation
to~$\Pi_n*\d_{-\lfloor nc\rfloor}$ provided by the centred Poisson distribution
$\widehat{\rm Po}(nv_c)$ is almost exactly the best that could be expected.
 
The proof is based on exploiting the equation
\begin{equation}\label{DF}
  {\mathbb E}({\mathcal A}_n h)(Z_{n})=0,
\end{equation}
where ${\mathcal A}_n$ denotes the infinitesimal generator of~$Z_{n}$:
$$
  ({\mathcal A}_n h)(i)\ :=\ 
   \sum_{j\in \mathbb Z\setminus \{0\}}n{{\lambda}_{j}}\Big(\frac{i}{n}\Big)\big[h(i+j)-h(i)\big],
\quad i \in {\mathbb Z},
$$
and where, here and subsequently, the quantity $Z_n$, when appearing without a time argument, 
is to be interpreted in such expressions as being a random variable having the
equilibrium distribution~$\Pi_n$.  The equation~\Ref{DF} is a manifestation of
Dynkin's formula, and it holds under rather mild conditions on~$h$: see
Hamza \& Klebaner~(1995).

Manipulations carried out in [SB] show that ${\mathcal A}_n h$ can be expressed in alternative form.

\begin{lemma}{\rm [SB, Lemma 1.1.]} \label{generator} 
Suppose that $\s^2(z) < \infty$ for all~$z\in\re$.
Then, for any function $h\colon \mathbb Z\to \mathbb R$  with bounded differences,  we have 
\eq\label{basic}
   ({\mathcal A}_nh)(i) 
      \ =\ \frac{n}{2}{\sigma}^2\Big(\frac{i}{n}\Big)\bigtriangledown{g_h}(i)
        +nF\Big(\frac{i}{n}\Big)g_h(i) + E_n(g_h,i), \non
\en
where $\bigtriangledown{f}(i):=f(i)-f(i-1)$ and $g_h(i) := \bigtriangledown{h}(i+1)$
and, for any $i\in\bZ$,
\eqa
   \lefteqn{E_n(g,i)}\non\\ 
  &:=&  -\frac{n}{2}F\Big(\frac{i}{n}\Big)\bigtriangledown{g}(i)
      +\sum_{j\geq 2}a_j(g,i) n\lambda_j\Big(\frac{i}{n}\Big) 
      -\sum_{j\geq 2}b_j(g,i) n{{\lambda}_{-j}}\Big(\frac{i}{n}\Big), 
    \label{En-def} \non
\ena 
with
\eqa
  a_j(g,i) &:=& -\binom{j}{2}\bigtriangledown {g}(i)  
     + \sum_{k=1}^{j-1} {k \bigtriangledown {g}(i+j-k)} \label{aj-bnd-1} \non \\ 
  &=&   \sum_{k=2}^j \binom{k}{2} \bigtriangledown^2{g}(i+j-k+1);\label{aj-bnd-2} \non \\
  b_j(g,i) &:=& \binom{j}{2}\bigtriangledown {g}(i)  
     - { \sum_{k=1}^{j-1} k \bigtriangledown {g}(i-j+k)} \non \\ 
  &=&   \sum_{k=2}^j \binom{k}{2} \bigtriangledown^2{g}(i-j+k). \non
\ena  
\end{lemma}

Since $F(c)=0$, we note that, for $i/n$ small, $\{-F'(c)\}^{-1}(\aa_n h)(i + \ncfl)$ 
is close to
\[
    \frac1{-F'(c)}\,\frac n2\,\s^2(c) \D g^*_h(i) - (i-\nvan)g^*_h(i)
		\Eq nv_c \D g^*_h(i) - (i-\nvan)g^*_h(i),
\]
for $g^*_h(i) := g_h(i+\ncfl)$, where $\nvan = nv_c - \nvfl$ denotes the fractional
part of~$nv_c$.  This is the Stein operator for the centred Poisson distribution
$\widehat{\rm Po}(nv_c)$ (R\"ollin,~2005), acting on the function~$g^*_h$. Combining this observation
with~\Ref{DF} and writing $Y_n=Z_n-\lfloor nc \rfloor$ yields
\eqa
  0 &=& \{-F'(c)\}^{-1}\ex\{(\aa_n h)(Y_n + \ncfl)\} \non\\
    &=& \ex\{nv_c\D g^*_h(Y_n) - (i-\nvan)g^*_h(Y_n)\} + \ex\{H(g^*_h,Y_n)\},
		\label{nearly-Stein}
\ena
say. If the error term $\ex\{H(g^*_h,Y_n)\}$ can be controlled,
then Stein's method leads easily to the approximation of
$\law(Y_n) = \Pi_n*\d_{-\ncfl}$ by $\widehat{\rm Po}(nv_c)$.
For the approximation of point probabilities, \Ref{nearly-Stein} needs to be 
analyzed for
functions~$g^*_h$ that are translates of the solutions to the Stein--Chen
equation corresponding to single point sets.

Carrying out this recipe, and examining the form of $H(g^*_h,Y_n)$, yields 
\begin{eqnarray}\nonumber
  \lefteqn{\sprz|(\Pi_{n}-\lfloor nc \rfloor)(r)-\widehat{{\rm Po}}(n v_c)(r)|} \\\non
  &&\Le \frac{1}{-F^{\prime}(c)} \sprz\; |{\mathbb E} R(n,r;Y_n)|
      +\sprz\;n v_c\; |\ex\{ \bigtriangledown^2{\tgr(Y_{n}+1)}\}|\\\nonumber
  &&\qquad\mbox{}+\sprz\;\widehat{{\rm Po}}(n v_c)\{r\}\cdot
     {\mathbb P}(Y_n< -\lfloor n v_c \rfloor)\\
  &&\ :=\ R_{n1}+R_{n2}+R_{n3}, \label{ek}
\end{eqnarray}\\
say, where 
\begin{eqnarray}
  R(n,r;Y_n)&:=& \frac{n}{2}\Bigl[{\sigma}^2\Bigl(\frac{Y_n+\lfloor nc \rfloor}{n}\Bigr)
     -\sigma^2(c)\Bigr]\bigtriangledown{\tgr(Y_n)} \nonumber \\
   &&\quad\mbox{} +n\Bigl[F\Bigl(\frac{Y_n+\lfloor nc \rfloor}{n}\Bigr)-F(c)
       -\frac{Y_n}{n}F^{\prime}(c)\Bigr]\tgr(Y_n) \nonumber \\ 
   &&\quad\mbox{} + F^{\prime}(c)\langle nv_c \rangle \tgr(Y_n) + E_n(\tgr,Y_n+\ncfl),
   \label{rt}
\ena
and the function~$\tgr$ is given by
\begin{equation}\label{defgt}
  \tgr(l)\ :=\ \left\{\begin{array}{ll}
       0,     &{\rm if}\;\;\; l < -\lfloor n v_c \rfloor \\
       g_{n v_c, \{r+\lfloor n v_c\rfloor\}}(l+\lfloor n v_c \rfloor), 
          &{\rm if} \;\;\; l\geq -\lfloor n v_c \rfloor. 
       \end{array} \right. 
\end{equation}
Here, for $A \subset \integ_+$, $g_{\m,A}$ denotes the solution to the 
Stein--Chen equation
\begin{equation}\label{SE}
   \charI_A(i)-{\rm Po}(\m)\{A\} \Eq \m\;g_{\m,A}(i+1)-i\;g_{\m,A}(i), \qquad i\geq 0.
\end{equation}
We further split the last term of~\Ref{rt} into
\[
   E_n(\tgr,Y_n+\ncfl) \Eq \sum_{l=1}^7 E_{nl}(\tgr,Y_n+\ncfl),
\]
with
{\allowdisplaybreaks
\eqa
  \label{En1}
    E_{n1}(\tgr,Y_n+\ncfl) &:=& -\frac{n}{2}\Bigl[F\Bigl(\frac{Y_n+\lfloor nc \rfloor}{n}\Bigr)
     -F(c)\Bigr]\btd\tgr(Y_n) ; \\
  E_{n2}(\tgr,Y_n+\ncfl)  &:=&   \sum_{j=2}^{\snl} \Bigl[\sum_{k=2}^{j}\binom{k}{2}\btd^2
       \tgr(Y_n+j-k+1)\Bigr]\,n \l_j(c); \label{En2} \\
   E_{n3}(\tgr,Y_n+\ncfl)  &:=&  \sum_{j=2}^{\snl} \Bigl[-\binom{j}{2}\btd\tgr(Y_n)
        +\sum_{k=1}^{j-1} k \btd\tgr(Y_n+j-k) \Bigr]\non\\
     &&\mbox{}\qquad n \Bigl\{\l_j\Bigl(\frac{Y_n+\lfloor nc \rfloor}{n}\Bigr) - \l_j(c)\Bigr\}; 
        \label{En3} \\
   E_{n4}(\tgr,Y_n+\ncfl)  &:=&  \sum_{j\ge \snu}\Bigl[-\binom{j}{2}\btd\tgr(Y_n)+\sum_{k=1}^{j-1} k
       \btd\tgr(Y_n+j-k)\Bigr]\non\\
      &&\mbox{}\qquad    n \l_j\Bigl(\frac{Y_n+\lfloor nc \rfloor}{n}\Bigr); \label{En4} \\
   E_{n5}(\tgr,Y_n+\ncfl)  &:=&  -\sum_{j=2}^{\snl}
     \Bigl[\sum_{k=2}^{j}\binom{k}{2}\btd^2\tgr(Y_n-j+k)\Bigr]\,n\l_{-j}(c); \label{En5}\\
   E_{n6}(\tgr,Y_n+\ncfl)  &:=&  -\sum_{j=2}^{\snl} \Bigl[\binom{j}{2}\btd\tgr(Y_n)
     -\sum_{k=1}^{j-1} k \btd\tgr(Y_n-j+k)\Bigr]\non\\
     &&\mbox{}\qquad n
       \Bigl\{\l_{-j}\Bigl(\frac{Y_n+\lfloor nc \rfloor}{n}\Bigr) - \l_{-j}(c)\Bigr\};  
			 \label{En6} \\
   E_{n7}(\tgr,Y_n+\ncfl)  &:=&  - \sum_{j\ge \snu}
	   \Bigl[\binom{j}{2}\btd\tgr(Y_n)-\sum_{k=1}^{j-1} k \btd\tgr(Y_n-j+k)\Bigr]\non\\
     &&\mbox{}\qquad n\lambda_{-j}\Bigl(\frac{Y_n+\lfloor nc \rfloor}{n}\Bigr).
         \label{En7}  
\ena
}
Our strategy for proving Theorem~\ref{limth} is now to show that each of the terms 
$R_{n1}$, $R_{n2}$ and $R_{n3}$ in~\Ref{ek} is 
of the desired order $O(n^{-(\a+1)/2}\sln)$; clearly, the treatment of~$R_{n1}$, which
involves all the detail of~$E_n(\tgr,Y_n+\ncfl)$, is to be the most laborious.

 \setcounter{equation}{0} 
 \setcounter{theorem}{0}  
\section[Point probabilities]{Differences of point probabilities}
As an essential step in proving Theorem~\ref{limth}, we need first to show that the differences
between the successive point probabilities of~$\Pi_n$ are suitably small.  
The bound that we achieve is of
order $O(n^{-1}\sqrt{\log n})$.  In order to prove this result, we begin with
two lemmas.  The first states that, for any $U\ge1$, the distribution of~$Z_n(U)$ 
has point probabilities
which are uniformly of order $O(n^{-1/2})$, if~$Z_n(0)$ is close enough to~$nc$.

\begin{lemma}\label{sdoi}
Under Assumptions A1--A4, for any $U\ge1$, there exists $C_{\ref{sdoi}}(U) < \infty$ such that 
$$
  \sup_{k\in \mathbb Z} \sup_{|i-nc|\le n\delta/2}\; {\mathbb P}(Z_{n}(U)=k \mid
     Z_{n}(0)=i)\ \leq\ C_{\ref{sdoi}}(U)n^{-1/2}.
$$
\end{lemma}

\begin{proof}
 Note that, for any integer valued random variable~$X$, 
\begin{eqnarray}
   \sup_{k\in\integ}\mathbb P(X=k) 
	 &=& \sup_{k\in\integ}\{\mathbb P(X\leq k)-\mathbb P(X+1\leq k)\}
      \non \\
   &\le& d_{TV}\{\mathcal{L}(X),\mathcal{L}(X)*\delta_1\}, \label{point-bnd}
\end{eqnarray} 
where $\mathcal{L}(X)$ denotes the distribution of $X$. Taking $X = Z_n(U)$ and applying
Lemma~\ref{lema4} completes the proof.
\end{proof}

The next lemma shows that the {\em differences\/} between successive point probabilities
of~$Z_n(U)$ are uniformly close, to order~$O(n^{-1}\sqrt{\log n})$,
for a large range of values of~$Z_n(0)$ and for a particular choice of~$U\ge1$.
This is the result
that we shall then be able to extend to the equilibrium distribution~$\Pi_n$.
For $\L^* := \sup_{|z-c| \le \d/2} \L(z)$, we set
\eq\label{U-def}
  U\ :=\ \max\{1,1/2\L^*\};\qquad \d'_1 \ :=\ \d e^{-U\|F'\|_\d}/4.
\en	

\begin{lemma}\label{sunu}
Under Assumptions A1--A4, and for $U$ and~$\d'_1$ defined above,
there exists $C_{\ref{sunu}}< \infty$ such that 
\begin{eqnarray*}
  &&\sup_{k\in \mathbb Z}\sup_{|i-nc|\le n\delta'_1} 
   |{\mathbb P}(Z_{n}(U)=k \mid Z_{n}(0)=i-1)-{\mathbb P}(Z_{n}(U)=k+1 \mid Z_{n}(0)=i)|\\
  &&\qquad \qquad \qquad \leq\ C_{\ref{sunu}}\ n^{-1}\sqrt{\log n}.
\end{eqnarray*} 
\end{lemma}

\begin{proof}
We compare the probability measures $\law((Z_n(u),\,0\le u\le U) \giv Z_n(0) = i-1)*\d_1$
and $\law((Z_n(u),\,0\le u\le U) \giv Z_n(0) = i)$ by examining the likelihood ratio of
the processes $Z_n\ui =_d \{Z_n \giv Z_n(0)=i-1\}$ and $Z_n\ut =_d \{Z_n \giv Z_n(0)=i\}$
 along paths with the same set of jumps $(j_l,\,l\ge1)$
occurring at the same times $(t_l,\,l\ge1)$.
$Z_n\ui$ starts from the state~$i-1$; we write
$z_l := n^{-1}\{i-1+\sum_{s=1}^l j_s\}$ for the value of $n^{-1}Z_n\ui$ at
time~$t_l$, $l\ge0$. $Z_n\ut$ starts from the state~$i$, and thus has the same
paths as $Z_n\ui + 1$.  Then
the likelihood ratio of the two processes along the first~$m$ steps of the path is given
by
\eqs
  S_m &:=& S_m(z_0,z_1\ldots,z_m; t_1,\ldots,t_m)\\
  &=& \prod_{l=0}^{m-1} \Blb \frac{\l_{j_{l+1}}(z_l + n^{-1})}{\l_{j_{l+1}}(z_l)}
   \exp\{-n(\L(z_l+n^{-1}) - \L(z_l))(t_{l+1} - t_l)\} \Brb \\
  &=&  \prod_{l=0}^{m-1} V_l.
\ens
Note that, since $|(1+x)(1+y) - 1| \le 3|x| + |y|$
if $|y| \le 2$, and since $|e^t-1| \le 2|t|$ in $t \le 1$, it follows that
\[
  |V_l-1| \Le  3\Bigl|\frac{\l_{j_{l+1}}(z_l + n^{-1})}{\l_{j_{l+1}}(z_l)} - 1 \Bigr|
   + 2n|\L(z_l+n^{-1}) - \L(z_l)|(t_{l+1} - t_l)
\]
provided that
\[
    n\{\L(z_l) - \L(z_l+n^{-1})\}(t_{l+1} - t_l) \Le 1.
\]

Now, if $|z-c| \le \d/2$ and $n^{-1} \le \d/2$, it follows from Assumptions~A2
and~A4 that
\[
  \Bigl|\frac{\l_{j}(z + n^{-1})}{\l_{j}(z)} - 1 \Bigr|
  \Le \frac{\| \l'_j\|_\d}{n\e\l_j(c)} \Le \frac{L_1}{n\e},
\]
and similarly that
\eq\label{L-ratio}
   \Bigl|\frac{\L(z + n^{-1})}{\L(z)} - 1 \Bigr|
  \Le \frac{\| \L'\|_\d}{n\e\L(c)} \Le \frac{L_1}{n\e}
\en
also.  Hence, for all $n\ge2/\d$, writing $e_{l+1} := n\L(z_l)(t_{l+1} - t_l)$, we have
\eq\label{Ul-bnd}
  |V_l-1| \Le \frac{L_1}{n\e}\{3 + 2e_{l+1}\},
\en
as long as 
\eq\label{1st-conds}
  |z_l-c| \le \d/2\quad \mbox{and either}\quad  \L(z_l) \le \L(z_l+n^{-1})
	     \quad \mbox{or}\quad e_{l+1} \le n\e/L_1.
\en

Now consider the random likelihood ratio process 
\[ 
   (S_m(n^{-1}Z_0,n^{-1}Z_1,\ldots,n^{-1}Z_m;\t_1,\ldots,\t_m),\, m\ge0), 
\]
where $(\t_l,\,l\ge0)$ denote the jump times of the process~$Z_n\ui$, 
and $Z_l := Z_n\ui(\t_l)$, $l\ge0$, the sequence of states 
that it visits; define also $E_{l} \ :=\ n\L(n^{-1}Z_{l-1})(\t_{l} - \t_{l-1})$. 
Then~$S := (S_m,\, m\ge0)$, is a martingale with mean~$1$ with respect
to the filtration $\gg_m := \s(Z_0,Z_1,\ldots,Z_m;\t_1,\ldots,\t_m)$, $m\ge0$.  We shall,
for technical reasons, work rather with another martingale~$\tS$, which 
typically agrees with~$S$
for a long time, but which satisfies the inequality
\eq\label{key-bnd}
  |\tS_{m+1} - \tS_m| \Le \frac{2L_1}{n\e}\,\{3 + 2E_{m+1}\}
\en
for all~$m\ge0$.  This we achieve by defining
$\s := \min\{\s_r,\,1\le r\le 3\}$, where
\eqa
  \s_1 &:=& \inf\{l\ge0 \colon\, n|\L(n^{-1}[Z_{l-1}+1]) - \L(n^{-1}Z_{l-1})|(\t_l - \t_{l-1}) > 1\},
                   \label{sig1-def}\\
  \s_2 &:=& \inf\{l\ge0 \colon\, S_l > 2\} \quad \mbox{and} \quad 
  \s_3 \ :=\ \inf\{l\ge0 \colon\, |n^{-1} Z_l - c| > \d/2\},\label{sig23-def}
\ena
and then setting
\[
   \tS_m \ :=\ S_{m\wedge\s} C_{m,\s_1},
\]
where
\[
    C_{m,l} \ :=\ \begin{cases}
       e/V_{l-1} &\mbox{if}\quad l \le \min\{m,\s_2,\s_3\}\ \mbox{and}\ 
			    \L(z_{l-1}) > \L(z_{l-1}+n^{-1});\\
		1 &\mbox{else}.		
       \end{cases}
\]
Note that the only effect of the factor~$C_{m,\s_1}$ is to multiply
$\tS$ by~$e$ instead of by~$V_{\s_1-1}$ at time~$\s_1$, if $\s_1 \le \min\{\s_2,\s_3\}$
and $\L(z_{\s_1-1}) > \L(z_{\s_1-1}+n^{-1})$. The value~$e$ is chosen so that the
martingale property is preserved; and the modification also ensures that~\Ref{key-bnd}
is still satisfied at time~$\s_1$, since $2(e-1)$ is no larger that $4L_1 E_{\s_1}/\{n\e\}$, 
because, at time $\s_1$,
\eqs
    1  &<&  n|\L(z_{\s_1-1}+n^{-1}) - \L(z_{\s_1-1})|(\t_{\s_1} - \t_{\s_1-1}) \\
     &=& E_{\s_1}\Bigl|\frac{\L(z_{\s_1-1}+n^{-1})}{\L(z_{\s_1-1})} - 1\Bigr|
     \Le   E_{\s_1}\,\frac{L_1}{n\e},
\ens
in view of~\Ref{L-ratio}.  

Now, from~\Ref{key-bnd}, and since also, by the strong Markov property, the
conditional distribution $\law(E_{l+1} \giv \gg_l)$ is the standard exponential
$\exp(1)$ distribution for each~$l$, the process~$\tS$ satisfies
the conditions of the variant of the bounded differences inequality for martingales
given in Barbour~(2008, Lemma~4.1), from which it follows that
\[
  \pr\left[|\tS_m - 1| > C\,\frac{L_1\sqrt{m\log m}}{n\e} \Giv Z_n(0) = i-1\right] 
     \Le 2\exp\{-3C\log m/928\}
\]
for any $m$ such that 
\[
     \sqrt{\frac m{\log m}} \ \ge\ 135 C / 236.
\]
In particular, recalling~\Ref{U-def}, for $m = m(n) := \lceil 2n\L^*U \rceil$,
we have
\eq\label{R-tilde-bnd}
    \pr\left[|\tS_{m(n)} - 1| > C\,\frac{L_1\sqrt{m(n)\log m(n)}}{n\e} \Giv Z_n(0) = i-1\right] 
            \Le 2n^{-3},
\en
if we take $C := 928$, as long as $n\ge e$ and
\eq\label{n-cond}
    \frac n{\log n} \ \ge\ 540^2.
\en

Now let $M_n(U) := \min\{l\colon\,\t_l > U\}$, and introduce the notation $\pr_s$
to denote $\pr[\cdot\giv Z_n(0)=s]$. Then
\eqa
   \lefteqn{\pr_{i-1}[\{M_n(U) > m(n)\} \cap \{\s_3 \ge M_n(U)\}] }\non\\
     &&\Eq \pr_{i-1}[\{\t_{m(n)} \le U\} \cap \{\s_3 \ge M_n(U)\}]\non\\
     &&\Le \Po(n\L^*U)\{(2n\L^*U,\infty)\} \Le \exp\{-n\L^*U/3\}, \label{chernoff}
\ena
by the Chernoff inequality (see Chung \& Lu (2006, Theorem~4)).  Hence, for~$U$ as defined
in~\Ref{U-def}, we have 
\eqa
  \lefteqn{|\pr_i[Z_n(U) = k+1 ] - \pr_{i-1}[Z_n(U) = k]|
      \Eq |\ex_{i-1}\{(S_{M_n(U)} - 1)I[Z_n(U) = k]\}| }\non\\
  &\le& |\ex_{i-1}\{(\tS_{M_n(U)} - 1)I[Z_n(U) = k] I[\tS_{M_n(U)} = S_{M_n(U)}]\}| 
     \   +\ \pr_{i-1}[\tS_{M_n(U)} \neq S_{M_n(U)}]\non\\
  &\le& |\ex_{i-1}\{(\tS_{m(n)} - 1)I[M_n(U) \le m(N)]I[Z_n(U) = k] I[\tS_{M_n(U)} = S_{M_n(U)}]
                      \}|\non\\
    &&\mbox{}  + \pr_{i-1}[\{\tS_{M_n(U)} \neq S_{M_n(U)}\} \cap \{M_n(U) \le m(N)\}]
                + \pr_{i-1}[M_n(U) > m(n)].\label{half-done}
\ena
Applying \Ref{R-tilde-bnd} and~\Ref{chernoff}, \Ref{half-done} now implies that, for
all~$n \ge e$ satisfying~\Ref{n-cond},
\eqa
   \lefteqn{|\pr_i[Z_n(U) = k+1] - \pr_{i-1}[Z_n(U) = k]|}\non\\
   &\le& \frac{CL_1}{n\e} \sqrt{m(n)\log m(n)}\, \pr_{i-1}[Z_n(U)=k] + 2n^{-3}
        + \pr_{i-1}[\cup_{l=0}^4 A_l] ,
       \label{nearly}
\ena 
where 
\eqs
   A_0 &:=& \Blb|\tS_{m(n)} - 1| > C\,\frac{L_1\sqrt{m(n)\log m(n)}}{n\e}\Brb; 
      \qquad   A_1 := \{\s_1 \le M_n(U)\};\\ 
   A_2 &:=& \{\s_2 \le M_n(U)\};\qquad A_3\ :=\ \{\s_3 < M_n(U)\};
   \qquad A_4\ :=\ \{M_n(U) > m(n)\}.
\ens
First, we note that $\pr_{i-1}[A_0] \le 2n^{-3}$, from~\Ref{R-tilde-bnd},
for all~$n$ such that~\Ref{n-cond} is satisfied.  Then,
from Lemma~\ref{lema3}, for all $|i - nc| \le n\d'_1$, as defined in~\Ref{U-def}, we have
\[
   \pr_{i-1}[A_3 ] 
      \Eq \pr_{i-1}\Bigl[\sup_{0\le u\le U} |n^{-1}Z_n(u) - c| > \d/2 \Bigr]
   \Le n^{-1}K_{U,\d/2};
\]
then
\[
   \pr_{i-1}[A_4 \cap A_3^c ] \Le \exp\{-n\L^*U/3\},
\]
by~\Ref{chernoff}; then, from~\Ref{sig1-def} and the definition of~$E_l$, and using~\Ref{L-ratio},
we have
\[
    A_1 \cap A_3^c \cap A_4^c \ \subset\ \bigcup_{l=1}^{m(n)} \Blb \frac{L_1 E_l}{n\e} > 1 \Brb ,
\]
so that
\[
   \pr_{i-1}[A_1 \cap A_3^c \cap A_4^c] \Le m(n) e^{-n\e/L_1}.
\]
Finally, we immediately have
$
   A_2 \cap A_1^c \cap A_3^c \cap A_4^c \ \subset\ A_0
$
for all~$n \ge \max\{3,2\L^*\}$ such that $n/\log n > 3(L_1/\e)^2$.  
Combining these bounds with~\Ref{nearly}, and noting also that,
from Lemma~\ref{sdoi},
\[
   \pr_{i-1}[Z_n(U) = k ] \Le C_{\ref{sdoi}}(U)/\sqrt n,
\]
for all~$|i-nc|\le n\delta'_1$,  the lemma is proved.   
\end{proof}

\begin{theorem}\label{supkpin}
Under Assumptions A1--A4, there exists a constant $C_{\ref{supkpin}}>0$ such that 
$$
   \sup_{k\in \mathbb Z} |\Pi_{n}(k)-\Pi_{n}(k+1)|\Le C_{\ref{supkpin}}\,n^{-1}\sqrt{\log n}.
$$ 
\end{theorem}

\begin{proof} Fix $U$ as in~\Ref{U-def}.
Since $\Pi_n$ is the equilibrium distribution of $Z_n$, it is in particular  true that
 \begin{eqnarray*}\nonumber
   \lefteqn{|\Pi_{n}(k)-\Pi_{n}(k+1)|
     \Eq |\sum_{i \in \mathbb Z}  \Pi_n(i){\mathbb P}_i(Z_{n}(U)=k )
       -\sum_{i \in \mathbb Z}  \Pi_n(i){\mathbb P}_{i}(Z_{n}(U)=k+1 )|}\\\nonumber
    &&\leq \ \sum_{i \in \mathbb Z}  \Pi_n(i-1)
      |{\mathbb P}_{i-1}(Z_{n}(U)=k )-{\mathbb P}_i(Z_{n}(U)=k+1)|\\
    &&\qquad\mbox{} +\sum_{i \in \mathbb Z} |\Pi_n(i-1)- \Pi_n(i)|\,
                         {\mathbb P}_i(Z_{n}(U)=k+1).\phantom{HHHHHHHHHHHH}
\end{eqnarray*} 
With $\delta'_1$ as in~\Ref{U-def}, note that one can write
 {\allowdisplaybreaks 
\begin{eqnarray*}\nonumber
   \lefteqn{\sum_{i \in \mathbb Z}  \Pi_n(i-1)\,|{\mathbb P}_{i-1}(Z_{n}(U)=k)
      -{\mathbb P}_i(Z_{n}(U)=k+1)|}\\
   &&\leq\ \Pi_n\Big\{|Z_n+1-nc|> n\delta'_1\Big\}
      +\sup_{|i-nc|\le n\delta'_1} |{\mathbb P}_{i-1}(Z_{n}(U)=k) 
          -{\mathbb P}_i(Z_{n}(U)=k+1)|,
\end{eqnarray*} }\\ 
 and that 
\begin{eqnarray*}\nonumber
 \lefteqn{\sum_{i \in \mathbb Z} |\Pi_n(i-1)- \Pi_n(i)|{\mathbb P}_i(Z_{n}(U)=k+1)}\\
   &&\le\ \Pi_n\Big\{|Z_n+1-nc| > n\delta'_1\Big\}+\Pi_n\Big\{|Z_n-nc| > n\delta'_1\Big\}\\
   &&\quad\mbox{} +\sup_{|i-nc|\le n\delta'_1} {\mathbb P}_i(Z_{n}(U)=k+1) \cdot 2
       d_{TV}\{\Pi_n,\Pi_n*\delta_1\}.
\end{eqnarray*} 
By applying the result of Corollary~\ref{prob-n1} three times, we obtain that
{\allowdisplaybreaks 
\begin{eqnarray}\nonumber
   \lefteqn{\sup_{k\in \mathbb Z}|\Pi_{n}(k)-\Pi_{n}(k+1)|}\\\nonumber
   &&\Le  O(n^{-1})  
     + \sup_{k\in \mathbb Z}\sup_{|i-nc|\le n\delta'_1} |{\mathbb P}_{i-1}(Z_{n}(U)=k)
      -{\mathbb P}_i(Z_{n}(U)=k+1)|\\\nonumber
   &&\qquad +\sup_{k\in \mathbb Z}\sup_{|i-nc|\le n\delta'_1} {\mathbb P}_i(Z_{n}(U)=k+1 ) 
     \cdot 2 d_{TV}\{\Pi_n,\Pi_n*\delta_1\}\\
   &&=:\ O(n^{-1})+\h_{1n}+\h_{2n}. \label{pk-pk1}
\end{eqnarray} }
The quantity $\h_{1n}$ is of order $O(n^{-1}\sln)$, in view of Lemma~\ref{sunu};
and Lemma~\ref{sdoi} and  Theorem~\ref{thdtvpi} together give the bound
\begin{eqnarray}\nonumber
  \h_{2n}
     &\le&  C_{\ref{sdoi}}(U)n^{-1/2}\cdot C_{\ref{thdtvpi}}\, n^{-1/2} \Eq O(n^{-1}).
		 \label{bonds2}
\end{eqnarray}   
This completes the proof of the theorem.
\end{proof}

 \setcounter{equation}{0} 
 \setcounter{theorem}{0}  
\section[Local limit theorem]{Local limit approximation for the equilibrium distribution}
The proof of Theorem~\ref{limth} consists of bounding
the quantities $R_{n1}$, $R_{n2}$ and~$R_{n3}$ of~\Ref{ek}, which all involve the
functions~$\tgr$ defined in~\Ref{defgt}.  For use in the subsequent argument, we collect
some of their properties.  We write $\|f\|_\infty := \sup_{i\in\integ}|f(i)|$,
$\|f\|_1 := \siz |f(i)|$.

\begin{lemma}\label{2.4}
We have the following estimates:
\eqs
  \brak1.&& \|\tgr\|_\infty \Le \|\D\tgr\|_\infty \Le 1/(nv_c);\\
  \brak2.&& \|\D \tgr(i)\|_1 \Le 2/(nv_c); \\
  \brak3.&& \|\D^2 \tgr(i)\|_1 \Le 4/(nv_c); \\
  \brak4.&& |(i-\langle nv_c\rangle) \tgr(i)| \Le h(i) + \Po(nv_c)\{r+\nvfl\}; \\
  \brak5.&& |(i - \langle nv_c\rangle) \D\tgr(i)| \Le h(i+1) + h(i) + 1/(nv_c),
\ens 
where, in parts 4 and~5, we have $h(i) \ge 0$ for all~$i$, and $\|h(i)\|_1 \le 3$.
\end{lemma}

\begin{proof}
For $i \le -\nvfl$, $\tgr(i) = 0$;  for $i > -\nvfl$, we have
$\tgr(i) = g_{\m,s}(j)$, where $j = i+\nvfl$, $\m = nv_c$ and $s = r+\nvfl$, 
and $g = g_{\m,s}$ satisfies the Stein--Chen equation
\eq\label{St-Poi}
   \m g(j+1) - jg(j) \Eq \bone_{\{s\}}(j) - \Po(\m)\{s\},\qquad j \ge 0.
\en
Parts 1 and~2 now follow from the proof of Lemma~1.1.1 of
Barbour, Holst \& Janson~(1992), in which it is shown that the function
$g_{\m,s}$ is negative and strictly decreasing in $\{1,2,\ldots,s\}$ and positive 
and strictly decreasing in $\{s+1,s+2,\ldots\}$, with $\D g_{\m,s}(s) \Le 1/(nv_c)$.
Part~3 is then immediate from part~2.

For part~4, using the notation above and~\Ref{St-Poi}, we have
\eqa
   \lefteqn{(i - \langle nv_c\rangle) \tgr(i) \Eq (j-\m)g_{\m,s}(j)} \non\\
   &&\Eq \m (g_{\m,s}(j+1) - g_{\m,s}(j)) - \bone_{\{s\}}(j) + \Po(\m)\{s\}. \label{i-factor}
\ena
This implies that
\[
  |(i - \langle nv_c\rangle) \tgr(i)| \Le \{\m|\D g(j)| + \bone_{\{s\}}(j)\} + \Po(\m)\{s\},
\]
which, with part~2, proves part~4.  It also follows immediately from \Ref{i-factor} that
\eqa
  |(i - \langle nv_c\rangle) \D\tgr(i)| \Le h(i+1) + h(i) + |\tgr(i+1)|,
\ena
for the same function~$h(i) := \{\m|\D g(j)| + \bone_{\{s\}}(j)\}$, and part~5 follows
on applying part~1. 
\end{proof}

As a result of these bounds, combined with Theorems~\ref{thdtvpi} and~\ref{supkpin},
we can establish two useful bounds on expectations of differences of the~$\tgr(Y_n+\cdot)$,
under the equilibrium distribution.

\begin{lemma}\label{new-one}
For any $r,l\in\integ$, we have
\eqs
  \brak1.&& \ex|\btd\tgr(Y_n+l)| \Le \frac{2C_{\ref{thdtvpi}}}{n^{3/2}v_c};\\
  \brak2.&& |\ex\{\btd^2\tgr(Y_n+l)\}| \Le \frac{2 C_{\ref{supkpin}}}{n^2v_c}\,\sln.
\ens
\end{lemma}

\begin{proof}
For the first part, it is immediate that
\[
    \ex|\btd\tgr(Y_n+l)| \Le \sup_{i'\in\integ}\Pi_n(i')\,\siz|\btd\tgr(i)|.
\]
By Lemma~\ref{2.4}\,(2) and~\Ref{point-bnd}, this is bounded in turn by
\[
    \dtv\{\Pi_n,\Pi_n*\d_1\}\,2/(nv_c),
\]
and part~1 follows from Theorem~\ref{thdtvpi}.  For the second part,
\begin{eqnarray}\label{term2insupk} \non
  \lefteqn{|{\mathbb E}\{ \btd^2 \tgr(Y_{n}+l)\}|} \\ \nonumber
  &&= \  \Bigl|\sum_{i\in \mathbb Z} \btd \tgr(i-\lfloor nc\rfloor+s)
                (\Pi_n(i+1)-\Pi_n(i))\Bigr|\\
  && \leq\ \Blb \sup_{i'\in \mathbb Z} \big|\Pi_{n}(i'-1)-\Pi_{n}(i')\big|\Brb\,  
      \sum_{i\in \mathbb Z} |\btd \tgr(i-\lfloor nc\rfloor)|\\\nonumber
  && \leq\  \sup_{i\in \mathbb Z} \big|\Pi_{n}(i-1)-\Pi_{n}(i)\big|\cdot 2(n v_c)^{-1},
\end{eqnarray}
where the last line uses Lemma~\ref{2.4}\,(2).  Part~2 of the lemma now follows from
Theorem~\ref{supkpin}.
\end{proof}

Bounding a further set of expectations that appear repeatedly in the estimates
first needs another, technical lemma.

\begin{lemma}\label{Lem2.5}
Let~$\m$ be any probability distribution on~$\integ$.  Suppose that $s,f$ and~$h$ are real
functions on~$\integ$ such that $\|f\|_\infty < \infty$,  $\|\D s\|_\infty < \infty$
and $\|h\|_1 < \infty$, which also satisfy the inequality 
\eq\label{sf-ineq}
  |s(i)f(i)| \Le |h(i)| + k, \qquad I_1 \le i < I_2, 
\en
 for some integers $I_1 < I_2$ and for some $k < \infty$.  Then
\eqs
  \Bigl|\sum_{i=I_1}^{I_2} \m_i s(i)\btd f(i)\Bigr| &\le& \|f\|_\ii\,\|\D s\|_\ii 
     + \|h\|_1 \sup_{I_1\le i < I_2}|\m_i-\m_{i+1}| + k\dtv(\m,\m*\d_1)\\
   &&\mbox{}\qquad  + |\m_{I_1}s(I_1)f(I_1-1)| + |\m_{I_2}s(I_2)f(I_2)|,
\ens
where $\|g\|_\ii := \sup_{I_1\le i < I_2}|g(i)|$. 
\end{lemma}

\begin{proof}
It is immediate that
\eqs
  \lefteqn{\Bigl|\sum_{i=I_1}^{I_2} \m_i s(i)\btd f(i)\Bigr|}\\
  && \Le 
       \Bigl|\siii\{\m_{i+1}s(i+1)-\m_is(i)\}f(i)\Bigr| + |\m_{I_1}s(I_1)f(I_1-1)|
           + |\m_{I_2}s(I_2)f(I_2)|\\
  &&\Le \Bigl|\siii\{\m_{i+1}-\m_i\}s(i)f(i)\Bigr| 
      + \Bigl|\siii\m_{i+1}\{s(i+1)-s(i)\}f(i)\Bigr|\\
   &&\mbox{}\qquad\qquad\qquad + |\m_{I_1}s(I_1)f(I_1-1)|
           + |\m_{I_2}s(I_2)f(I_2)|.
\ens
Clearly, the second term is bounded by $\|f\|_\ii\,\|\D s\|_\ii$.  For the
first term, in view of~\Ref{sf-ineq}, we have at most
\[
  \siii \{|\m_{i+1}-\m_i|\,|h(i)|\} +  k\siii|\m_{i+1}-\m_i|,
\]
which is easily bounded by $\|h\|_1 \sup_{I_1 \le i < I_2}|\m_i-\m_{i+1}| + k\dtv(\m,\m*\d_1)$,
in view of~\Ref{point-bnd}.

Note that the argument also goes through for $I_1=-\infty$ and $I_2=\infty$, in which
case the final two elements in the bound disappear.
\end{proof}

This lemma is combined with Lemma~\ref{2.4}\,(4) and~(5) to give the next
corollary, which is used as an ingredient in many of the estimates to be made.  

\begin{corollary}\label{sf-appns}
Suppose that $|s(i)| \le |i - \ncfl|$ for all~$|i| \le n\d$.  Then, for any $0 < \d'\le \d$ and
all $l\in\integ$ such that $|l| \le n(\d-\d')$, we have
\eqs
  \lefteqn{1.\quad |\ex\{s(Y_n+l)\btd\tgr(Y_n+l)I[|Y_n| \le n\d']\}| }\\
   &&\Le \frac1{nv_c}\,\sup_{|i|\le n\d}|\D s(i)| 
      + \frac{3C_{\ref{sunu}}}n \sqrt{\log n} + \frac{C_{\ref{thdtvpi}}}{2n\sqrt{v_c}}
      + 2(C_{\{\ref{2ineq},1\}} + C_{\{\ref{2ineq},2\}}/\d')/(nv_c);\\
  \lefteqn{2.\quad |\ex\{s(Y_n+l)\btd^2\tgr(Y_n+l)I[|Y_n| \le n\d']\}|}\\
   &&  \Le \frac2{nv_c}\,\sup_{|i|\le n\d}|\D s(i)|
      + \frac{6C_{\ref{sunu}}}n \sqrt{\log n} + \frac{C_{\ref{thdtvpi}}}{n^{3/2}v_c}
      + 4(C_{\{\ref{2ineq},1\}} + C_{\{\ref{2ineq},2\}}/\d')/(nv_c).
\ens
\end{corollary}

\begin{proof}
We take $\Pi_n*\d_{-l}$ for~$\m$ and either $\tgr$ or $\btd\tgr$ for~$f$ in Lemma~\ref{Lem2.5},
noting that parts 4 and~5 of Lemma~\ref{2.4} give the appropriate counterparts
of~\Ref{sf-ineq}.  The first three elements appearing in the bound given by Lemma~\ref{Lem2.5}
are in turn bounded by using Lemma~\ref{2.4}\,(1), Lemma~\ref{sunu} and Theorem~\ref{thdtvpi}.
The last two are bounded by Lemma~\ref{2.4}\,(1) and Theorem~\ref{exubj}.
\end{proof}

\medskip
\nin{\it Proof of Theorem~\ref{limth}}

\medskip
We are now in a position to undertake the proof of
Theorem~\ref{limth}, for which we need to bound the terms $R_{1n}$, $R_{2n}$ and $R_{3n}$ in~\Ref{ek}
to order $O(n^{-(\a+1)/2}\sln)$. First, we show that $R_{3n}$ is as
small as $O(n^{-3/2})$. This is because, from Barbour \& Jensen~(1989, Remark to Lemma 2.1), 
if $X\sim\Po(\m)$, then
$$
  \sup_{k\in \mathbb Z}\; \mathbb P (X=k)\leq \frac{1}{2\sqrt{\mu}}.
$$ 
Hence, and from Corollary~\ref{prob-n1}, it follows easily that
\begin{equation}\label{R3-bnd}
  R_{3n} \Eq \sup_{k\in \mathbb Z} \;\widehat{{\rm Po}}(n v_c)\{k\}
    \cdot{\mathbb P}(Y_n< -\lfloor n v_c \rfloor)
    \Eq O\Big(\frac{1}{n\sqrt{n}}\Big).
\end{equation}
For the quantity~$R_{2n}$ in~\Ref{ek}, we just use Lemma~\ref{new-one}\,(2) to give
\eq\label{R2-bnd}
   R_{2n} \ :=\ nv_c \sup_{r\in\integ}|{\mathbb E}\{ \bigtriangledown^2 \tgr(Y_{n}+1)\}|
     \Le 2C_{\ref{supkpin}}\,n^{-1}\sqrt{\log n}.
\en
It thus remains to bound~$R_{1n}$.  To do so, we consider in turn
the expectations of the quantities appearing in \Ref{rt} and in \Ref{En1}--\Ref{En5}.

Beginning with the elements of $\ex R(n,r;Y_n)$, we first have
\eq\label{Rnr-1}
  \ex\Blb\frac{n}{2}\Bigl[{\sigma}^2\Bigl(\frac{Y_n+\lfloor nc \rfloor}{n}\Bigr)
     -\sigma^2(c)\Bigr]\btd{\tgr(Y_n)}\Brb,
\en
which is of the form considered in Corollary~\ref{sf-appns}\,(1), with $l=0$ and
\eqs
   s(i) &:=&  \frac{n}{2}\Bigl[{\sigma}^2\Bigl(\frac{i+\lfloor nc \rfloor}{n}\Bigr)
     -\sigma^2(c)\Bigr].
\ens
For $|i| \le n\d/2$ and $n \ge 2/\d$, we have
\[
   |s(i)| \Le \half\,|i-\langle nv_c\rangle|\,\|(\s^2)'\|_\d \quad\mbox{and}
   \quad |s(i)-s(i-1)| \le \half\,\|(\s^2)'\|_\d,
\]
whereas, for $|i| > n\d/2$, we have the simple bound
\[ 
   |s(i)| \Le \frac{n}{2}\Bigl[\s^2(c) + \sjz j^2c_j(1 + n^{-1}|i|)\Bigr],
\]
using Assumption~A2. By Theorem~\ref{exubj} and Corollary~\ref{prob-n1}, it follows
that the latter element contributes at most $O(n^{-1})$ to $|\ex R(n,r;Y_n)|$; for the
former, Corollary~\ref{sf-appns} gives a bound of order $O(n^{-1}\sln)$.

For the next term,
\[
   \ex \Blb n\Bigl[F\Bigl(\frac{Y_n+\lfloor nc \rfloor}{n}\Bigr)-F(c)
       -\frac{Y_n}{n}F^{\prime}(c)\Bigr]\tgr(Y_n)\Brb,
\]
$|\tgr(Y_n)|$ is bounded by $1/(nv_c)$, using Lemma~\ref{2.4}\,(1).
The contribution from the part $|Y_n| \le n\d$ is thus easily bounded by 
\[
   \bigl[\|F''\|_\d\, n^{-2}\,\ex \{Y_n^2 I[|Y_n \le n\d]\} + \|F'\|_\d n^{-1}\bigr]/v_c,
\]
and $\ex \{Y_n^2 I[|Y_n \le n\d]\}= O(n)$
by Theorem~\ref{exubj}, so that the whole contribution is of order $O(n^{-1})$.
If $|Y_n| > n\d$, Assumption~A2 and Theorem~\ref{exubj} guarantee a contribution
of the same order.  The third term immediately yields 
\[
   \ex|F^{\prime}(c)\langle nv_c \rangle \tgr(Y_n)| \Le |F^{\prime}(c)|/(nv_c),
\]
again of order $O(n^{-1})$.
All of these elements are of order $O(n^{-1}\sln)$, at least as small as the order
$O(n^{-(1+\a)/2}\sln)$ stated in the theorem, and it thus
remains to bound $|\ex\{E_{nl}(\tgr,Y_n+\ncfl)\}|$ for $1\le l\le 7$.

For the term arising from~\Ref{En1}, we have
\[
   \ex\Bigl\{\frac{n}{2}\Bigl[F\Bigl(\frac{Y_n+\lfloor nc \rfloor}{n}\Bigr)
     -F(c)\Bigr]\btd\tgr(Y_n)\Bigr\},
\]
which is of the form considered in Corollary~\ref{sf-appns}\,(1), with $l=0$ and
\eqs
   s(i) &:=&  \frac{n}{2}\Bigl[F\Bigl(\frac{i+\lfloor nc \rfloor}{n}\Bigr)
     -F(c)\Bigr],
\ens 
and can be treated very much as was~\Ref{Rnr-1}, yielding a bound of the
same order.  For that arising from~\Ref{En2}, 
\[
  \ex \Blb\sum_{j=2}^{\snl} \Bigl[\sum_{k=2}^{j}\binom{k}{2}\btd^2
       \tgr(Y_n+j-k+1)\Bigr]\,n \l_j(c)\Brb,
\]
we can use Lemma~\ref{new-one}\,(2) 
to bound the expectations $\ex\btd^2\tgr(Y_n+j-k+1)$, giving a contribution of at most 
\[
   \sum_{j=2}^{\snl}\frac16\, j^3c_j\,n \frac{2C_{\ref{supkpin}}}{n^2v_c}\,\sqrt{\log n}
      \Le \frac{C_{\ref{supkpin}}\,s_\a}{3v_c}\,n^{-(1+\a)/2}\sln,
\]
where we have also used Assumption~A2. 

The next term is from~\Ref{En3}, and is more complicated.  For its summands, we write
\eqa
   \lefteqn{\Bigl[-\binom{j}{2}\btd\tgr(Y_n)+\sum_{k=1}^{j-1} k \btd\tgr(Y_n+j-k)\Bigr]\,
       n \Blb\l_j\Bigl(\frac{Y_n+\lfloor nc \rfloor}{n}\Bigr) - \l_j(c)\Brb}\non \\
  && \Eq -\binom{j}{2}\btd\tgr(Y_n)\,n 
      \Blb\l_j\Bigl(\frac{Y_n+\lfloor nc \rfloor}{n}\Bigr) - \l_j(c)\Brb \non\\
  &&\qquad\mbox{} + \sum_{k=1}^{j-1} k \btd\tgr(Y_n+j-k)\,n
      \Blb\l_j\Bigl(\frac{Y_n+\lfloor nc \rfloor+j-k}{n}\Bigr) - \l_j(c)\Brb \non\\
  &&\qquad\mbox{} + \sum_{k=1}^{j-1} k \btd\tgr(Y_n+j-k)\,n
      \Blb\l_j\Bigl(\frac{Y_n+\lfloor nc \rfloor}{n}\Bigr) 
         - \l_j\Bigl(\frac{Y_n+\lfloor nc \rfloor+j-k}{n}\Bigr)\Brb \non\\
  && \Eq E_{n3}\ui(Y_n,j) +  E_{n3}\ut(Y_n,j) +  E_{n3}\uh(Y_n,j), \label{En3l}
\ena
say.  The term $E_{n3}\ui(Y_n,j)$ is of the form considered in Corollary~\ref{sf-appns}\,(1), 
with $l=0$ and
\eqs
   s(i) &:=&  - n\binom{j}{2} 
      \Blb\l_j\Bigl(\frac{i+\lfloor nc \rfloor}{n}\Bigr) - \l_j(c)\Brb.
\ens
For $|i| \le n\d/2$,
\[
   |s(i)| \Le \binom{j}{2}\,|i-\langle nv_c\rangle|\,\|\l_j'\|_\d \quad\mbox{and}
   \quad |s(i)-s(i-1)| \Le \binom{j}{2}\,\|\l_j'\|_\d,
\]
whereas, for $|i| > n\d/2$, we have the direct bound
\[
   |s(i)| \Le nc_j\binom{j}{2} (2 + n^{-1}|i|),
\]
using Assumption~A2.  From Corollary~\ref{sf-appns} and Assumption~A4, the contribution
from the first part is of order 
\eq\label{En31-bnd}
   O\Bl c_j\binom{j}{2} n^{-1}\sln \Br;
\en
the second part is also at most of this order, in view of Theorem~\ref{exubj}, 
Corollary~\ref{prob-n1} and~Lemma~\ref{2.4}\,(2).  Adding over $j \le \snl$,
this gives a total contribution to the quantity $|\ex\{E_{n3}(\tgr,Y_n+\ncfl)\}|$ of order
$O(n^{-1}\sln)$.

For $E_{n3}\ut(Y_n,j)$, we now have a sum of terms of the form considered in 
Corollary~\ref{sf-appns}\,(1), with $l=j-k$ and
\eqs
   s(i) &:=& nk \Blb\l_j\Bigl(\frac{i+\lfloor nc \rfloor}{n}\Bigr) - \l_j(c)\Brb.
\ens
Supposing~$n$ to be large enough that $\sqrt n \le n\d/2$, we have
\[
   |s(i)| \Le k\|\l_j'\|_\d\,|i-\langle nv_c\rangle| \quad\mbox{and}
   \quad |s(i)-s(i-1)| \Le k\|\l_j'\|_\d
\]
for $|i| \le n\d/2$, whereas, for $|i| > n\d/2$, we have the bound
\[
   |s(i)| \Le nkc_j (2 + n^{-1}|i|).
\]
Arguing very much as for~\Ref{En31-bnd}, it thus follows that the total contribution 
to the quantity $|\ex\{E_{n3}(\tgr,Y_n+\ncfl)\}|$ is again of order $O(n^{-1}\sln)$.

Finally, for $E_{n3}\ut(Y_n,j)$, we again have a sum of terms.  We first note that
\eqs
   \Bigl|\l_j\Bigl(\frac{i+\lfloor nc \rfloor}{n}\Bigr) 
         - \l_j\Bigl(\frac{i+\lfloor nc \rfloor+j-k}{n}\Bigr)\Bigr|
     &\le& n^{-1}|j-k|\|\l_j'\|_\d
\ens
for $|i| \le n\d/2$, and this leads to a contribution to~$|\ex E_{n3}\ut(Y_n,j)|$ of
at most
\eq\label{En33-1}
   \sum_{k=1}^{j-1} k(j-k)\|\l_j'\|_\d/(nv_c) \Le L_1j^3c_j/(6nv_c),
\en
in view of Lemma~\ref{2.4}\,(1).  For $|i| > n\d/2$, there is the bound
\[
   |E_{n3}\ut(i,j)| \Le \sum_{k=1}^{j-1} \frac k{v_c}c_j\{2 + n^{-1}(2|i| + j-k)\},
\]
giving
\eq\label{En33-2}
   |\ex\{E_{n3}\ut(Y_n,j) I[|Y_n| > n\d/2]\}| \Le j^2c_j\{{\textstyle\frac76}\, \pr[|Y_n| > n\d/2]
      + n^{-1}\ex(|Y_n|I[|Y_n| > n\d/2])\},
\en
because $j \le \snl$.  Adding \Ref{En33-1} and~\Ref{En33-2} over $j \le \snl$ gives
a total contribution to $|\ex E_{n3}\ut(Y_n,j)|$ of order $O(n^{-(1+\a)/2})$, because
of Assumption~A2.

The term from~\Ref{En4} is much easier.  For $|i| \le n\d$, we have the bound
\[
   \l_j\Bigl(\frac{i+\lfloor nc \rfloor}{n}\Bigr) \Le c_j(1 + \d),
\]
by Assumption~A2,
and $\ex|\btd \tgr(Y_n+l)| \le 2(C_{\ref{thdtvpi}}/v_c) n^{-3/2}$ for any~$l$, by 
Lemma~\ref{new-one}. Hence
\eqs
  \lefteqn{\ex\Bigl|\Bigl[-\binom{j}{2}\btd\tgr(Y_n)+\sum_{k=1}^{j-1} k
       \btd\tgr(Y_n+j-k)\Bigr]\,n \l_j\Bigl(\frac{Y_n+\lfloor nc \rfloor}{n}\Bigr)
          I[|Y_n| \le n\d]\Bigr|}\\
  && \Le  j^2 c_j (1+\d) \frac{2C_{\ref{thdtvpi}}}{v_c}\, n^{-1/2},
   \phantom{HHHHHHHHHHHHHHHHHHHHHHHH}
\ens
and summing over $j \ge \snu$ gives a total contribution to~\Ref{En4} of at most
\[
     s_\a \frac{2C_{\ref{thdtvpi}}}{v_c}\,(1+\d) n^{-(1+\a)/2},
\]
in view of Assumption~A2.
For $|Y_n| > n\d$, the $j$-contribution is bounded by 
\[
    j^2 c_j \|\btd\tgr\|_\infty \ex\{(n + |Y_n|)I[|Y_n| \ge n\d]\} 
       \Le 2 j^2 c_j (C_{\{\ref{2ineq},1\}} + C_{\{\ref{2ineq},2\}})/(nv_c),
\]
in view of Lemma~\ref{2.4}\,(1) and Theorem~\ref{exubj}, and summing over $j \ge \snu$ gives
a contribution of order $O(n^{-1-\a/2})$.  Hence the complete contribution from~\Ref{En4}
is of order $O(n^{-(1+\a)/2})$.

The remaining terms \Ref{En5}--\Ref{En7} are treated in exactly the same way as those
in \Ref{En2}--\Ref{En4}.  In all, the largest order of any of the terms in \Ref{En1}--\Ref{En7}
is of order $O(n^{-(1+\a)/2}\sln)$, and since the other terms were of order $O(n^{-1}\sln)$, 
Theorem~\ref{limth} is proved. \hfill$\qed$

 \setcounter{equation}{0} 
 \setcounter{theorem}{0}
\setcounter{section}{4} 
\section*{Appendix} \label{appendix}
The following results from~[SB] are used in the proofs.  
    
\begin{theorem}{\rm [SB, Theorem 2.1.]} 
\label{exubj} Under Assumptions A1--A4, for all~$n$  large enough, the process $Z_n$ has an
   equilibrium distribution $\Pi_n$, and
\begin{equation} \label{2ineq}
  \begin{split}
    &{\mathbb E}\{|n^{-1}Z_n-c| \cdot \charI(|n^{-1}Z_n-c| > \delta)\}
              \ \le\ C_{\{\ref{2ineq},1\}}n^{-1};\\
    &{\mathbb E}\{(n^{-1}Z_n-c)^2 \cdot \charI(|n^{-1}Z_n-c|\le \delta)\}
              \ \le\ C_{\{\ref{2ineq},2\}}n^{-1},
  \end{split}
\end{equation}
for $\delta$ as in Assumption~A3 and constants $C_{\{\ref{2ineq},1\}}$ and~$C_{\{\ref{2ineq},2\}}$;
as before, in such expressions, $Z_n$ is used to denote a random variable having the
equilibrium distribution~$\Pi_n$.
\end{theorem}

\begin{corollary} {\rm [SB, Corollary 2.5.]}  \label{prob-n1} 
Under Assumptions A1--A4, for any fixed $\d'$ such that $0 < \d' \le \d$,
there exists $C_{\ref{prob-n1}}(\d') < \infty$ such that
$$
   {\mathbb P}[|n^{-1}Z_n-c|> \d'] \Le  C_{\ref{prob-n1}}(\d')n^{-1}.
$$
\end{corollary}

\begin{lemma}{\rm [SB, Lemma 3.1.]}\label{lema3}
Under Assumptions A1--A4, for any $U > 0$ and $0 < \h \le \d$, there exists a constant 
$K_{U,\h}<\infty$ such that
$$
   \mathbb P[\sup_{t\in [0,U]}|Z_n(t)-nc|> n\h  \mid Z_n(0) = i]\ \leq\  n^{-1}K_{U,\h},
$$ 
uniformly in $|i - nc| \le (n\h/2) \exp\{-\|F'\|_{\d}U\}$.
\end{lemma}

\begin{theorem} {\rm [SB, Theorem 3.2.]} \label{thdtvpi} 
Under Assumptions A1--A4, there exists a constant $C_{\ref{thdtvpi}}>0$ such that
$$
   d_{TV}\{\Pi_{n},\Pi_{n}*\delta_1\} \Le C_{\ref{thdtvpi}} n^{-1/2},
$$ 
where $\Pi_{n}*\delta_1$ denotes the unit translate of $\Pi_{n}$.
\end{theorem} 

Finally, we shall use the following result, which was used in~[SB] to prove the
previous theorem; see, for example, (3.7) in the proof of [SB, Theorem 3.2].

\begin{lemma}\label{lema4}
Under Assumptions A1--A4,  for any $U\ge1$, there exists a constant $K_U<\infty$ such that
$$
   d_{TV}\{\law(Z_n(U) \giv Z_n(0)=i),\law(Z_n(U) \giv Z_n(0)=i)*\delta_1\} \Le K_U n^{-1/2},
$$ 
uniformly in $|i - nc| \le n\delta/2$.
\end{lemma}


\begin{thebibliography}{00}

\bibitem{AB80} 
{\sc A.\ D.\ Barbour} (2008)\; 
Coupling a branching process to an infinite dimensional epidemic process.
Preprint

\bibitem{BHJ92} 
{\sc A.\ D.\ Barbour, L.\ Holst \& S.\ Janson} (1992)\; 
{\it Poisson Approximation\/}.  Oxford Univ. Press.

\bibitem{BX99} 
{\sc A.\ D.\ Barbour \& J.\ L.\ Jensen} (1989)\; 
Poisson perturbations. {\it ESAIM,  P\&S~{\bf 3}}, 131--150.

\bibitem{LC75} 
{\sc L.\ H.\ Y.\ Chen} (1975)\; 
Poisson approximation for dependent trials. {\it Ann.\ Probab.~{\bf 3}}, 534--545.

\bibitem{cl06}
{\sc F.\ Chung \& L.\ Lu} (2006)\;
Concentration inequalities and martingale inequalities: a survey.
{\it Internet Math.\/}~{\bf 3}, 79--127. 

\bibitem{HK95} 
{\sc K.\ Hamza \& F.\ C.\ Klebaner} (1995)\; 
Conditions for integrability of Markov chains. 
{\it J.Appl.\ Prob.~{\bf 32}}, 541--547.

{\sc \bibitem{TK81} 
T.\ G.\ Kurtz} (1981)\; 
{\it Approximation of population processes.\/} CBMS-NSF Regional Conf. Series
in Appl.\ Math.~{\bf 36},  SIAM, Philadelphia.
 
\bibitem{idi05} 
{\sc A.\ R\"ollin} (2005)\; 
Approximation of sums of conditionally independent random variables by the translated Poisson distribution. 
{\it Bernoulli~{\bf 11}},  1115--1128.

\bibitem{SB08} 
{\sc S.\ Socoll \& A.\ D.\ Barbour} (2008)\; 
Translated Poisson approximation to equilibrium distributions of Markov
population processes. Preprint.

\end{thebibliography}
\end{document}